\documentclass[12pt]{article}
\pdfoutput=1
\usepackage[utf8]{inputenc}
\usepackage{fullpage}
\usepackage[small,bf]{caption}

\usepackage{float}
\usepackage{geometry}
\usepackage{verbatim}
\usepackage{graphicx}
\usepackage{enumerate, enumitem}
\usepackage{booktabs}
\usepackage[
colorlinks=true,
citecolor=blue,
linkcolor=magenta,
urlcolor=cyan]{hyperref}
\usepackage{url}
\usepackage{xcolor}
\usepackage{tikz}
\usepackage{xparse}
\usepackage{amsmath,amsfonts}
\usepackage{bm}
\newcommand{\eg}{{\it e.g.}}
\newcommand{\ie}{{\it i.e.}}

\newcommand{\BEQ}{\begin{equation}}
\newcommand{\EEQ}{\end{equation}}
\newcommand{\BEAS}{\begin{eqnarray*}}
\newcommand{\EEAS}{\end{eqnarray*}}
\newcommand{\ones}{\mathbf 1}
\newcommand{\reals}{{\mbox{\bf R}}}

\newcommand{\symm}{{\mbox{\bf S}}}  % symmetric matrices
\newcommand{\diag}{\mathop{\bf diag}}

\newcommand{\intr}{\mathop{\bf int}}

\usepackage{relsize} % Relsize package needed

\usepackage{listings}

\definecolor{codegreen}{rgb}{0,0.6,0}
\definecolor{codegray}{rgb}{0.5,0.5,0.5}
\definecolor{codepurple}{rgb}{0.58,0,0.82}
\definecolor{backcolour}{rgb}{0.95,0.95,0.92}

\lstdefinestyle{mystyle}{
    backgroundcolor=\color{backcolour},   
    commentstyle=\color{codegreen},
    keywordstyle=\color{magenta},
    numberstyle=\tiny\color{codegray},
    stringstyle=\color{codepurple},
    basicstyle=\ttfamily\small,
    breakatwhitespace=false,         
    breaklines=true,                 
    captionpos=b,                    
    keepspaces=true,                 
    numbers=left,                    
    numbersep=5pt,                  
    showspaces=false,                
    showstringspaces=false,
    showtabs=false,                  
    tabsize=2
}

\definecolor{seagreen}{rgb}{0.18, 0.55, 0.34}
\definecolor{mediumviolet-red}{rgb}{0.78, 0.08, 0.52}
\definecolor{khaki}{rgb}{0.94, 0.9, 0.55}

\lstdefinelanguage{mypython}
{
	keywords=[1]{from, import, assert, not, print},
	keywordstyle=[1]{\color{mediumviolet-red}},
	keywords=[2]{surecr, torch, cp, lo, pl},
	keywordstyle=[2]{\color{seagreen}},
	numbers=none,
	upquote=true,
	showstringspaces=false,
	basicstyle=\ttfamily,
	columns=fullflexible,
	keepspaces=true,
	emph={True,False,as,def,return,float,class,match,switch,len},
	emphstyle={\color{seagreen}},
	frame=trBL,
	belowskip=1em,
	aboveskip=1em,
	captionpos=b
}

\usetikzlibrary{math}

\usepackage{textcomp}

\newcommand{\data}{\theta}

\begin{document}

\title{Differentiating Through a Quadratic Cone Program}
\author{
	Quill Healey\thanks{healeyquill@gmail.com}
	\and
	Parth Nobel\thanks{\{ptnobel, boyd\}@stanford.edu, Department of Electrical Engineering, Stanford University}
	\and
	Stephen Boyd\footnotemark[2]
}
\date{\today}

\maketitle

\begin{abstract}
Quadratic cone programs are rapidly becoming the standard canonical 
form for convex optimization problems.
In this paper we address the question of differentiating the solution
map for such problems, generalizing previous work for linear
cone programs.  We follow a similar path, using the implicit function
theorem applied to the optimality conditions for a
homogenous primal-dual embedding.
Along with our proof of differentiability, we present methods for
efficiently evaluating the derivative operator and its adjoint at
a vector. Additionally, we present an open-source implementation of
these methods, named \texttt{diffqcp}, that can execute on CPUs and GPUs.
GPU-compatibility is already of consequence as it enables convex
optimization solvers to be integrated into neural networks with reduced
data movement, but we go a step further demonstrating
that \texttt{diffqcp}'s performance on GPUs surpasses the performance
of its CPU-based counterpart for larger quadratic cone programs.
\end{abstract}

\clearpage

\tableofcontents

\clearpage

\section{Introduction}

A \textit{quadratic cone program} (QCP) is an optimization problem which minimizes
a convex quadratic function over the intersection of a subspace and a convex cone.
Quadratic cone programming is the generalization of both quadratic programming and (linear)
cone programming, which date to the 1950s~\cite{frank1956} and
1990s~\cite[Chapter 4]{nesterov1994}, respectively.
Specifically, a \textit{quadratic program} (QP) is a QCP whose cone is restricted to the
product of $\{0\}$, $\reals$, and $\reals_+$, while a (linear) \textit{cone program}
is a QCP restricted to having a linear objective.

Quadratic programs, despite their limited modeling power,
have been studied extensively as they arise ubiquitously across many disciplines---from classical
engineering contexts to finance.
Cone programming, on the other hand, has been studied for its generality---all convex
optimization problems can be equivalently written as a cone program. 
Along with their rich theory, signficant development has gone into specialized solvers
for both quadratic programs~\cite{OSQP,proxqp} and cone programs~\cite{ECOS,SCS}.
Moreover, domain specific languages, such as CVXPY \cite{CVXPY} and CVXR \cite{CVXR},
have been designed to enable easy modeling with both classes of programs~\cite{CVX}.

Perturbation and sensitivity analysis has also been thoroughly developed
for QPs and cone programs.
Clasically, this analysis centered on the Lagrange multipliers~\cite{rockafellar1970convex,rockafellar1998variational,robinson1981sensitivity}.
In recent years, \textit{differentiable optimization}---the
derivative of the solution map between an optimization problem's parameters and its solution---has been developed.
The gradients of the solution map of a quadratic program (with respect to the problem data) were derived in~\cite{amos2017optnet}
by exploiting the problem structure. Subsequently,~\cite{Agrawal2019DifferentiatingTA}
proposed a technique for differentiating the solution map of a cone program using a more general
approach based on the implicit function theorem.
Differentiable optimization has found applications across energy systems~\cite{degleris2024gradient},
statistics~\cite{SURECR,randalo}, control~\cite{agrawal2020learning, Barratt2019FittingAK},
and in neural networks~\cite{cvxpylayers2019,negiar2023}.

In recent years, specialized solvers for QCPs have been developed and have
demonstrated significant speedups on problems previously solved via cone programs
or quadratic programs~\cite{SCS3,Goulart2024ClarabelAI}.
As a result, QCPs are emerging as a practical alternative to cone programs for many
convex optimization problems. Further, there has been success at GPU-accelerating
these QCP solvers~\cite{chen2024cuclarabel}.
However, the theory of differentiating the solution map of QCPs has remained undeveloped.

\subsection{Our contribution}
Closely following \cite{Agrawal2019DifferentiatingTA}, we derive conditions for when
the derivative, and its adjoint, of the primal-dual solution map to a QCP with
respect to the QCP's parameters exists.
We then present an extension to \cite{Agrawal2019DifferentiatingTA} to evaluate
Jacobian-vector and vector-Jacobian products
with these derivatives via projection onto cones and sparse linear system solves.
We then describe our GPU-accelerated Python implementation of this
method in \S\ref{sec:implementation}, which forms the derivative of the solution
map as an abstract linear operator.
Additionally, in \S\ref{sec:appendix-A} we present a unified reference of cone projection
operators and their derivatives.
Notably, this reference and our implementation includes the power cone,
which has previously been neglected in the differentiable optimization literature.

%\subsection{Automatic differentiation (AD)}
%
%\subsection{GPU Acceleration}

\section{Solution map and its derivative}

Following \cite{Agrawal2019DifferentiatingTA},
we consider the mapping from the numerical
data defining the primal and dual problems of a QCP to its solutions.
This \textit{solution map} is in general set-valued,
but in neighborhoods where it is single-valued
it is an implicit function of the problem data.
In the sequel, we present a system of equations % namely N(z, D)
that implicitly define the solution map of a QCP
when it is single-valued.
Applying the implicit function theorem to this
system, we obtain regularity conditions
on the problem data that
guarantee when the solution map is single-valued
and its derivative exists.
Finally, we provide an expression for the derivative
at points where these conditions are satisfied.

\subsection{QCPs and implicit functions}
The primal and dual problems for a (convex) QCP are
% \begin{equation}
%     \begin{array}{lll}
%         \text{(P)} &\text{minimize} \; & (1/2)x^T P x + q^T x  \\
%         &\text{subject to} & Ax + s = b  \\
%         & & s \in \mathcal{K},
%     \end{array}
%     \qquad
%     \begin{array}{lll}
%          \text{(D)}  &\text{maximize} \; & -(1/2)x^T P x -b^T y  \\
%         &\text{subject to} & Px + A^T y = -q \\
%         & & y \in \mathcal{K}^*,
%     \end{array}\label{eq:conic-pair}
% \end{equation}
\begin{equation}
    \begin{array}{lll}
        \text{(P)} &\text{minimize} \; & \tfrac{1}{2}x^T P x + q^T x  \\
        &\text{subject to} & Ax + s = b  \\
        & & s \in \mathcal{K},
    \end{array}
    \qquad
    \begin{array}{lll}
         \text{(D)}  &\text{maximize} \; & -\tfrac{1}{2}x^T P x -b^T y  \\
        &\text{subject to} & Px + A^T y = -q \\
        & & y \in \mathcal{K}^*,
    \end{array}\label{eq:conic-pair}
\end{equation}
where $x \in\reals^n$ is the \textit{primal} variable,
$y \in \reals^m$ is the \textit{dual} variable,
and $s \in \reals^m$ is the primal \textit{slack} variable.
We assume that $\mathcal K \subseteq \reals^m$ is a nonempty, closed,
convex cone with \textit{dual cone} $\mathcal K^*$.
The \textit{problem data} are $P\in \symm_+^{n}$, $A \in \reals^{m \times n}$,
$q \in \reals^n$, and $b \in \reals^m$.
(The convex cone can also be problem data, but for our purposes
we fix $\mathcal K$.) To simplify the subsequent discussion,
we define the set
\[
\Theta = \left\{
\data = (P, A, q, b)
\, \middle | \,
(P, A, q, b) \in \symm^{n} \times \reals^{m \times n} \times \reals^{n} \times \reals^{m}
\right\}.
\]
That is, $\data$ is the concatenation of problem data
(relaxed to allow $P \not \succeq 0$)---a
change from \cite{Agrawal2019DifferentiatingTA} which embeds the problem data
into a skew-symmetric matrix.
% \qoh{I added the ``relaxed''
% bit as I think we need to acknowledge
% why we say P is PSD when we define the problem data
% but then only require it to be symmetric in the concatenation.
% }

\paragraph{Optimality conditions.}
The optimality conditions for~\eqref{eq:conic-pair} are
\begin{equation}
Ax + s = b, \quad Px + A^T y = -q,\quad s \in \mathcal{K}, \quad  y \in \mathcal{K}^*, \quad s^T y = 0.
\label{eq:kkt-conds}
\end{equation}
Note that $s^T y$ is the duality gap, \ie, at any point that satisfies the first four
equalities and inclusions,  $s^Ty = \hat p - \hat d$ where $\hat p$ is the primal
objective at $(x, s)$ and $\hat d$ is the dual objective value at $(x, y)$.
Also note that~\eqref{eq:kkt-conds} is an implicit system that defines
the solution map to a QCP when it is single-valued.

\paragraph{Homogenous embedding.}
By applying $s^T y = \hat p - \hat d$,
a solution that satisfies~\eqref{eq:kkt-conds}
is equivalent to a solution of the following nonlinear systems of
equations
\begin{equation}
  \begin{bmatrix}
      -q \\ s \\ 0
  \end{bmatrix}
  =
  \begin{bmatrix}
      P & A^T  \\ -A & 0  \\ -q^T & -b^T
  \end{bmatrix}
  \begin{bmatrix}
      x \\ y
  \end{bmatrix} +
  \begin{bmatrix}
      0 \\ b \\ -x^T P x
  \end{bmatrix},
  \quad (x, s, y) \in \mathbf{R}^{n} \times \mathcal{K} \times \mathcal{K}^*.
  \label{eq:nonlin-system}
\end{equation}
However, because this system is not guaranteed to be feasible
(\eg, when the problem is primal or dual infeasible),
we instead consider the \textit{homogeneous embedding}
(as defined in~\cite{Goulart2024ClarabelAI})
\begin{equation}
  \begin{aligned}
    &\begin{bmatrix}
        0 \\ s \\ \kappa
    \end{bmatrix}
    =
    \begin{bmatrix}
        Px + A^T y + \tau q \\
        -Ax + \tau b \\
        -(1/\tau)x^T P x - q^T x - b^T y
    \end{bmatrix}, \\
    &(x, s, y, \tau, \kappa) \in \mathbf{R}^{n} \times \mathcal{K} \times \mathcal{K}^* \times \mathbf{R}_+ \times \mathbf{R}_+,
    \quad \tau + \kappa > 0,
  \end{aligned}
  \label{eq:nonlin-system2}
\end{equation}
where $\tau$ and $\kappa$ are new real-valued variables.
Unlike~\eqref{eq:nonlin-system},
this embedding is guaranteed to be (asymptotically) feasible even when~\eqref{eq:conic-pair}
is primal or dual infeasible.

Applying a change of variable with $N = n + m + 1$, the sets
\[
K = \reals^n \times \mathcal K^* \times \reals_+,
\quad
K^* = \{0\}^n \times \mathcal K \times \reals_+,
\]
the variables,
\[
  u = (x, y, \tau) \in \reals^N,
  \quad
  v = (0, s, \kappa) \in \reals^N
\]
and the functions $Q_1: \reals^N \to \reals^n$, $Q_2: \reals^N \to \reals^m$, and
$Q_3: \reals^N \to \reals$ defined as
\[
Q_1(u) = Px + A^Ty + \tau q, \quad
Q_2(u) = -Ax + \tau b, \quad
Q_3(u) = -(1/\tau)x^T P x - q^Tx - b^T y,
\]
we simplify the sequel by writing~\eqref{eq:nonlin-system2} as
\begin{equation}\label{eq:embedding}
  Q(u) = v, \quad u \in K, \quad v \in K^*, \quad u_N + v_N > 0.
\end{equation}
Here $Q: \reals^N \to \reals^N$ is defined as $Q(u) = (Q_1(u), Q_2(u), Q_3(u))$.
Lastly, we define a solution to \eqref{eq:embedding} as a complementary solution if $u_N v_N = 0$.

\subsection{Solution map}
% The solution map of a family of parameterized optimization problems is a
% (possibly multi-valued) mapping between the parameter values and
% the primal-dual solution to the problem.
For given problem data, the corresponding QCP~\eqref{eq:conic-pair} may have no solution,
a unique solution, or multiple solutions.
For the remainder of this paper, we assume it has a unique solution.
We define the solution map $S: \Theta \to \reals^{n+2m}$
of a family of parameterized optimization problems as the function
mapping $\theta$ to vectors $(x, y, s)$ that satisfy~\eqref{eq:kkt-conds}.
Similar to \cite{Agrawal2019DifferentiatingTA},
we express this function as composition of functions.
Unlike in \cite{Agrawal2019DifferentiatingTA},
we only have two functions in our composition: $S = \phi \circ s$, where
\begin{itemize}
\item $s: \Theta \to \reals^N$ maps the problem data to a
complementary solution of the homogeneous embedding and
\item $\phi: \reals^N \to \reals^{n + 2m}$ maps a complementary solution
to a solution of the primal-dual pair.
\end{itemize}

At a point $\data$ where $S$ is differentiable, the derivative of the
solution map is
\[
D {S}(\data) = D\phi \left( s(\data) \right) Ds( \data),
\]
by the chain rule.
In the remainder of this section we develop an expression for $DS(\data)$
by following the approach taken in~\cite{Agrawal2019DifferentiatingTA}:
\begin{itemize}
  \item We pose the problem of finding a (complementary) solution to the homogeneous
 embedding~\eqref{eq:embedding} as finding a root of a (differentiable) map,
 a function of both an input to the embedding and the primal-dual pair's problem data
 $\data$.
%  (This root is the solution to the implicit
%  equations that define the solution map to a QCP.)
  \item We consider the differentiability of this map and collect its derivatives with
 respect to both the embedding input and problem data.
  \item Using the implicit function theorem,
  we find $Ds(\data)$ in terms of these derivatives.
While $Ds(\data)$ will require the evaluation of $s$ at a point $\data$
and we never find an expression for $s$ directly, in practice we can supply such a
point, $s(\data)$, by using a
(convex) quadratic conic optimization numerical solver.
\end{itemize}

\subsection{Other machinery}

This subsection closely follows \cite{Busseti2018SolutionRA}.

\paragraph{The conic complementarity set.} The \textit{conic complementarity set} is defined as
\[
\mathcal{C} = \left\{ (u, v) \in K \times K^* 
\, \middle | \, u^T v = 0 \right\}.
\]
Let $\Pi$ and $\Pi^\circ$ be the projections onto
the cone $K$ and its polar cone $K^\circ = -K^*$,
respectively.
Note the functional equality (a form of the Moreau decomposition) \(\Pi^\circ = I - \Pi\),
where $I$ is the identity operator.

\paragraph{Minty's parameterization of the complementarity set.}
Let $M: \mathbf{R}^{N} \to \mathcal{C}$ be the Minty parameterization of $\mathcal{C}$, defined as
\[
M(z) = (\Pi z, - \Pi^\circ z),
\]
with inverse $M^{-1}: \mathcal{C} \to \mathbf{R}^{N}$ given by
\[
M^{-1}(u, v) = u - v.
\]

Unlike in~\cite{Busseti2018SolutionRA},
\begin{equation}
  -\Pi^\circ z = Q (\Pi z), \quad z_N \not = 0
  \label{eq:minty-restrictions}
\end{equation}
are not equivalent to the homogeneous embedded conditions~\eqref{eq:embedding}.
While $z$ satisfying \eqref{eq:minty-restrictions} implies that $u, v = M(z)$ satisfy
\eqref{eq:embedding},
there exists a non-complementary solution $(u, v)$ to \eqref{eq:embedding} such that
$z = M^{-1}(u, v)$ does not satisfy \eqref{eq:minty-restrictions}.
However, \eqref{eq:minty-restrictions} is equivalent to the
KKT conditions~\eqref{eq:kkt-conds}.

\paragraph{Residual map.}
The residual map $\mathcal{R}: \mathbf{R}^{N} \to \mathbf{R}^{N}$ is defined as
\begin{equation}
\mathcal{R}(z) = Q \left(\Pi z \right) + \Pi^\circ z = Q \left(\Pi z \right) - \Pi z + z.
\label{eq:residual-map}
\end{equation}
The map $\mathcal{R}$ is positive homogeneous and differentiable almost everywhere.

\paragraph{Normalized residual map.}
The \textit{normalized residual map}
$\mathcal{N}: \left\{ z \in \mathbf{R}^{N} \, \middle | \, z_N \not = 0 \right\} \to \mathbf{R}^{N}$
is defined as
\begin{equation}
  \mathcal{N}(z) = \mathcal{R}(z/\left| z_N \right|) = \mathcal{R}(z) / \left| z_N \right|.
  \label{eq:normalized-residual}
\end{equation}
By~\eqref{eq:minty-restrictions},
if $z \in \mathbf{R}^{N}$ is a solution to the conic pair~\eqref{eq:conic-pair}
then $\mathcal{N}(z) = 0$.
Conversely, if $\mathcal{N}(z) = 0$, then $z$ is a solution to \eqref{eq:conic-pair}.

\paragraph{Data dependence.}
Throughout this note we have fixed the data defining the primal-dual pair~\eqref{eq:conic-pair},
and consequently have not made explicit the dependence of $\mathcal R, \mathcal N, Q$
on $\data$. As we consider the derivative of these functions with respect to $\data$,
we will update our notation writing
\[
Q(u, \data), \qquad
\mathcal{R}(z, \data), 
\quad \text{and} \quad \mathcal{N}(z, \data).
\]

\subsection{Derivatives}

\paragraph{Normalized residual map with respect to the data.}
The normalized residual map is an affine function of the problem data, $\data$.
Therefore, it is differentiable with 
\[
D_{\data}\mathcal{N}(z, \data) [\widetilde{\data}] =
\frac{1}{\left| z_N \right|}D_{\data} Q(\Pi z, \data) [\widetilde{\data}]
\quad \text{and} \quad
D_{\data}\mathcal{N}(z, \data)^T [w]
= \frac{1}{\left| z_N \right|} D_{\data} Q(\Pi z, \data)^T [w],
\]
where
\[
D_{\data}Q(u, \data)[\widetilde{\data}] =
\begin{bmatrix}
  \widetilde{P} x  + \widetilde{A}^T y + \tau \widetilde{q} \\
  -\widetilde{A} x + \tau \widetilde{b} \\
  (-1/\tau)x^T \widetilde{P} x - \widetilde{q}^T x - \widetilde{b}^T y
\end{bmatrix}
\]
and $D_{\data}Q(u, \data)^T [w] = (\widetilde{P}, \widetilde{A}, \widetilde{b},
\widetilde{q}, \widetilde{b})$ where
\begin{align*}
  \widetilde{P} &= 1/2\left( w_{1:n} x^T + x w_{1:n}^T \right) - (w_N/\tau) x x^T,
  &\quad \widetilde{A} &= y w_{1:n}^T - w_{n+1:n+m} x^T, \\
  \widetilde{q} &= \tau w_{1:n} - w_N x,
  &\quad \widetilde{b} &= \tau w_{n+1:n+m} - w_N y,
\end{align*}
% \[
% D_{\data}Q(u, \data)^T [w] =
% \left(
% \begin{array}{c}
%   1/2\left( w_{1:n} x^T + x w_{1:n}^T \right) - (w_N / \tau) x x^T, \\[0.5em]
%   y w_{1:n}^T - w_{n+1:n+m} x^T, \\[0.5em]
%   \tau w_{1:n} - w_N x, \\[0.5em]
%   \tau w_{n+1:n+m} - w_N y
% \end{array}
% \right)
% \left(
for $w \in \reals^N$ and $\widetilde{\data} \in \Theta$.
\paragraph{Normalized residual map with respect to the variables.}
The normalized residual map is differentiable at $z$
if $z_N \neq 0$ and $\Pi$ is differentiable at $z$.
When $z$ is a solution of the primal-dual pair~\eqref{eq:conic-pair},
\[
D_{z}\mathcal{N}(z, \data) = \frac{1}{z_N} \left( D_{z}Q(\Pi z, \data) D\Pi(z) - D \Pi(z) + I \right),
\]
where the Jacobian of the nonlinear, homogeneous map $Q$ with respect to the embedding input is
\[
D_{u}Q(u, \data) = \begin{bmatrix}
  P & A^T & q \\
  -A & 0 & b \\
  (-2/\tau) x^T P - q^T & - b^T & (1/\tau^2) x^T P x
\end{bmatrix}.
\]

\paragraph{Implicit function theorem applied to $\mathcal N$.}
If $z$ is a solution of the primal-dual pair~\eqref{eq:conic-pair}
and $\Pi$ is differentiable at $z$,
then $\mathcal{N}$ is differentiable at $z$, $N(z, \data) = 0$, and $z_N > 0$.
Now suppose that $D_{z} \mathcal{N}(z, \data)$ is invertible.
The implicit function theorem \cite{rockafellar1998variational} guarantees that there
exists a neighborhood $V \subseteq \Theta$ of $\data$ on which the
solution $z = s(\data)$ of $\mathcal{N}(z, \data)$ is unique.
Furthermore, $s$ is differentiable on $V$,
$\mathcal{N}(s(\data), \data) = 0$ for all $\data \in V$, and
\[
Ds(\data) = - \left( D_z \mathcal{N}(z, \data)\right)^{-1} D_{\data}\mathcal{N}(z, \data).
\]

\paragraph{Solution construction.}
To construct a solution $(x, y, s)$ of the primal-dual pair~\eqref{eq:conic-pair} from
a complementary solution $z$ of the homogeneous embedding,
we use the function $\phi$ given in~\cite{Agrawal2019DifferentiatingTA}.
With  $\phi: \mathbf{R}^{N} \to \mathbf{R}^{n + 2m}$ given by
\[
\phi(z) = (z_{1:n}, \Pi_{\mathcal{K}^*}(z_{n+1:n+m}), \Pi_{\mathcal{K}^*}(z_{n+1:n+m}) - z_{n+1:n+m})/z_{N}.
\]
If $\Pi_{\mathcal{K}^*}$ is differentiable  at $z_{n+1:n+m}$, then $\phi$ is also differentiable and
\[
D\phi(z) = \begin{bmatrix}
    I & 0 & -x \\
    0 & D\Pi_{\mathcal{K}^*}(z_{n+1:n+m}) & -y \\
    0 & D\Pi_{\mathcal{K}^*}(z_{n+1:n+m}) - I & -s
\end{bmatrix}.
\]

\section{Implementation}\label{sec:implementation}

\subsection{Computing the Jacobian-vector product}\label{sec:jacobian-vector-product}

Applying the derivative $D \mathcal{S} (\data)$ to a perturbation
$d \data = (dP, dA, dq, db) \in \Theta$ corresponds to evaluating
\begin{align*}
    (dx, dy, ds) = D \mathcal{S} (\data) [d \data]
       &= D \phi (s(\data)) D s(\data) [d \data] \\
       &= D \phi (z) \left(- \left( D_z \mathcal{N}(z,
	\data)\right)^{-1} \right)
	D_{\data} \mathcal{N}(z, \data) [d \data].
\end{align*}
Given a solution $(x, y, s)$ to~\eqref{eq:conic-pair},
we construct a root of the normalized residual map as
$z = s(\data) = M^{-1}(u, v) = u - v$ where $u = (x, y, 1)$ and $v = (0, s, 0)$.

We now work from right to left.
First, we compute $\Pi z$ and form
$d_\data \mathcal{N} =  D_{\data} \mathcal{N}(z, \data) [d \data]$.
Second, we compute
\[
dz = - F^{-1} d_{\data}\mathcal{N},
\]
where
\[
F = D_z \mathcal{N}(z, \data) = 1/z_N \left(D_z Q( \Pi z, \data) D\Pi(z) - D\Pi (z) + I \right).
\]
Since it is impractical to form or factor $F$ as a dense matrix in some applications
(\textit{e.g.}, when $F$ is large), we use LSMR \cite{lsmr} to solve
\[
\underset{dz}{\text{minimize}} \quad \left\lVert Fdz + d_\data \mathcal N \right\rVert_{2}^{2},
\]
which only requires multiplication with $F$ and $F^T$.
Finally, we compute the solution perturbations as
\[
\begin{bmatrix}
    dx \\ dy \\ ds
\end{bmatrix}
=
\begin{bmatrix}
  dz_{1:n} - (dz_N) x \\
    D\Pi_{\mathcal{K}^*}(z_{n+1:n+m})[dz_{n+1:n+m}] - (dz_{N})y \\
    D\Pi_{\mathcal{K}^*}(z_{n+1:n+m})[dz_{n+1:n+m}] - (dz)_{n+1:n+m} - (dz_N)s
\end{bmatrix}.
\]

\subsection{Computing the vector-Jacobian product}
The adjoint of the derivative applied to a perturbation $(dx, dy, ds)$ is
\begin{align*}
  d\data = (dP, dA, dq, db) &= D \mathcal{S}(\data)^T [(dx, dy, ds)] \\
  &= Ds(\data)^T D\phi(s(\data))^T [(dx, dy, ds)] \\
  &= D_{\data}\mathcal{N} (z, \data)^T \left(- \left(D_{z} \mathcal{N} (z, \data)^T \right)^{-1} \right)D\phi(z)^T [(dx, dy, ds)],
\end{align*}
letting $z = s(\data)$ as in \S\ref{sec:jacobian-vector-product}.
Working right to left, first, we evaluate
\[
dz = D\phi(z)^T [(dx, dy, ds)] =
\begin{bmatrix}
  dx \\
  D\Pi_{\mathcal{K}^{*}}(z_{n + 1 : n+m})[dy + ds] - ds \\
  -x^T dx - y^T dy - s^T ds
\end{bmatrix}.
\]
Second, we evaluate $\Pi z$ and form $d_{\data} \mathcal{N} = - F^{-T} dz$ using LSMR.
Finally the problem data perturbation $d\data = (dP, dA, dq, db)$ is given by
% \begin{align*}
%   dP &= \frac{1}{2}\left( d \mathcal{N}_{n} \widetilde{z}_n^T + \widetilde{z}_n d \mathcal{N}_{n}^T \right) - (d \mathcal{N}_{N} / \widetilde{z}_N ) \widetilde{z}_n \widetilde{z}_n^T \\
%   dA &= \widetilde{z}_m d \mathcal{N}_{n}^T - d \mathcal{N}_{m} \widetilde{z}_n^T \\
%   dq &= \widetilde{z}_{N} d \mathcal{N}_{n} - d \mathcal{N}_{N} \widetilde{z}_{n} \\
%   db &= \widetilde{z}_{N} d \mathcal{N}_{m} - d \mathcal{N}_{N} \widetilde{z}_{m},
% \end{align*}
\begin{align*}
  dP &= \frac{1}{2}\left( d_{\data} \mathcal{N}_{1:n} \left(\Pi z \right)_{1:n}^T + \left(\Pi z \right)_{1:n} d_{\data} \mathcal{N}_{1:n}^T \right)
  - (d_{\data} \mathcal{N}_{N} / \left(\Pi z \right)_{N} ) \left(\Pi z \right)_{1:n} \left(\Pi z \right)_{1:n}^T, \\
  dA &= \left(\Pi z \right)_{n+1:n+m} d_{\data} \mathcal{N}_{1:n}^T - d_{\data} \mathcal{N}_{n+1:n+m} \left(\Pi z \right)_{1:n}^T, \\
  dq &= \left(\Pi z \right)_{N} d_{\data} \mathcal{N}_{1:n} - d_{\data} \mathcal{N}_{N} \left(\Pi z \right)_{1:n}, \\
  db &= \left(\Pi z \right)_{N} d_{\data} \mathcal{N}_{n+1:n+m} - d_{\data} \mathcal{N}_{N} \left(\Pi z \right)_{n+1:n+m}.
\end{align*}
% where
% \[
% d \mathcal{N}_{n} = d_{\data} \mathcal{N}_{1:n}, \qquad
% d \mathcal{N}_{m} = d_{\data} \mathcal{N}_{n+1:n+m},
% \quad \text{and} \quad
% d \mathcal{N}_{N} = d_{\data} \mathcal{N}_{N},
% \]
% and
% \[
% \widetilde{z}_n = \left(\Pi z \right)_{1:n}, \qquad
% \widetilde{z}_m = \left(\Pi z \right)_{n+1:n+m},
% \quad \text{and} \quad
% \widetilde{z}_N = \left(\Pi z \right)_{N}.
% \]
However, we do not form $dP$ and $dA$ exactly as formulated.
Instead we only compute their nonzero (or, more precisely, non-explicit-zero) entries as dictated
by the sparsity patterns of $P$ and $A$ respectively.

\subsection{Hardware accelerated Python implementation}
We have developed an open-source JAX~\cite{jax2018github} library
(also making significant use of the packages Equinox~\cite{kidger2021equinox}
and Lineax~\cite{lineax2023}), \texttt{diffqcp},
which implements these algorithms and is available at
\url{https://github.com/cvxgrp/diffqcp}.
Our implementation supports any QCP whose cone can be expressed
as the Cartesian product of the zero cone, the positive orthant,
second-order cones, and positive semidefinite cones. (Support for
exponential cones, power cones, and their duals is in development.)

% We use JAX~\cite{jax2018github} and libraries
% Equinox~\cite{kidger2021equinox} and Lineax~\cite{lineax2023}
% to enable
% computing the JVPs and VJPs on a GPU, allowing
% \texttt{diffqcp} to be integrated into GPU workflows,
% such as neural network training, without host-to-device
% data transfers.
\paragraph{Data movement.}
Host-to-device transfers have been a long-standing limitation of CVXPYlayers,
a Python library for constructing differentiable convex optimization
layers in PyTorch~\cite{paszke2019pytorch}, JAX, and
TensorFlow~\cite{tensorflow2015-whitepaper} using CVXPY.
Since CVXPYlayers only supports computing the derivative (and its
adjoint) of the solution map of a conic program on the CPU,
using this library to embed a differentiable convex optimization layer
in a neural network requires tranferring any data on the device
to the host during the forward or backward pass.
Such transfers can be expensive, so having to perform them on both
the forward and backward passes during every training iteration
can make this embedding prohibitive.
Being a JAX library, \texttt{diffqcp} can compute JVPs and VJPs on a
GPU, allowing our software to be integrated into GPU workflows, such
as neural network training, without these significant host-to-device
data transfers.

\paragraph{Performance.}
\texttt{diffqcp} relies on JAX to enable its high performance.
The JAX library uses Python as a “metaprogramming language” to build
performant and just-in-time compiled XLA programs.
Moreover we rely on the JAX transformation \texttt{vmap},
to simplify writing SIMD computations.
\texttt{diffqcp} makes extensive use of this transformation to ``batch''
projections onto a family of cones with the same dimensionality.
This batching is especially advantageous when computing JVPs and VJPs on a GPU,
as it enables the execution of many independent computations in parallel,
thereby maximizing processor occupancy and overall throughput.

\subsection{Example}\label{sec:example}
To test our implementation, we applied gradient descent to a loss function of the form
\[
	\|x - x^\star\|_2^2 + \|r - r^\star\|_2^2 +
	\|s - s^\star\|_2^2 + \|y - y^\star\|_2^2.
\]
where each $(x, r), s, y$ are the optimal primal, slack, and dual solutions of
\begin{equation}
      \begin{array}{ll}
	      \text{minimize} & r^T P r + q^T \begin{bmatrix}x \\ r \end{bmatrix} \\
	\text{subject to} & \begin{bmatrix}
		C & D  \\
		E & 0 \\
		f^T & 0^T \\
		\end{bmatrix} \begin{bmatrix} x \\ r \end{bmatrix} +s = b \\
	      & s \in \{0\}^m \times \reals_+^n \times \{0\}
      \end{array}
      \label{eq:prob-problem}
\end{equation}
	with $D, E, P$ are diagonal, and $q, f, b$ are vectors and where $(x^\star, r^\star), s^\star, y^\star$ are
the primal, slack, and dual solutions of 
\begin{equation*}
      \begin{array}{ll}
        \text{minimize} & r^T I r \\
	\text{subject to} & \begin{bmatrix}
		C^\star & -I  \\
		-I & 0   \\
		\ones^T & 0^T \\
		\end{bmatrix} \begin{bmatrix} x \\ r \end{bmatrix} +s = \begin{bmatrix}
	d^\star \\
	0 \\
	1
	\end{bmatrix} \\
	      & s \in \{0\}^m \times \reals_+^n \times \{0\}
      \end{array}
\end{equation*}
for a randomly selected $C^\star , d^\star$.

\paragraph{Results.} We take $m=2000$ and $n=1000$.
CuClarabel and \texttt{diffqcp} on an Intel Xeon E5-2670 CPU
and an NVIDIA TITAN Xp GPU took $44.20$ seconds per iteration.
As a control, we canonicalized the objective with SOCs and ran the gradient descent with Clarabel and \texttt{diffcp}, which took
$96.86$ seconds per iteration on the Intel Xeon E5-2670 CPU.
The improved modeling capacity and GPU-acceleration enabled a $2.19 \times$ speedup.

\clearpage
\appendix

\section{Cones, projections, and derivatives}\label{sec:appendix-A}

\paragraph{Zero cone.}
The zero cone $\{0\}$ has dual cone $\reals$, projection operation $\Pi_{\{0\}}(z) = 0$, and $D\Pi(z)[dz] = 0$.

\paragraph{Nonnegative cone.}
The nonnegative cone $\{x \mid x \geq 0\}$ has dual cone $\{x \mid x \geq 0\}$,
projection operator 
\[
\Pi_{\{x \mid x \geq 0\}}(z) = \begin{cases}
z & z \geq 0 \\
0 & z < 0,
\end{cases}
\]
and 
\[
D\Pi_{\{x \mid x \geq 0\}}(z)[dz] = \begin{cases}
dz & z > 0 \\
0 & z < 0.
\end{cases}
\]

\paragraph{Second-order cone.}
The second order cone $\mathcal K_{\mathrm{soc}} = \{(t, u) \in \reals \times \reals^n: \|u\|_2 \leq t\}$
has dual cone $\mathcal K_{\mathrm{soc}}$, projection operator
\[
\Pi_{\mathrm{soc}}((t, u)) = \begin{cases}
0 & \|u\|_2 \leq -t \\
(t, u) & \|u\|_2 \leq t \\
(1/2)(1 + t/\|u\|_2)(\|u\|_2, u) & \|u\|_2 \ge |t|
\end{cases}
\]
and
\[
D\Pi_{\mathrm{soc}}((t, u))[(dt, du)] = \begin{cases}
  0 & \|u\|_2 < -t \\
  (dt, du) & \|u\|_2 < t \\
  \frac{1}{2 \|u\|_2} \begin{bmatrix}
  \|u\|_2 dt + u^T du \\
  u dt + (t + \|u\|_2)du - (t / \|u\|_2^2) (u^T du) u
\end{bmatrix} & \|u\|_2 > |t|.
\end{cases}
\]

\paragraph{Positive semidefinite cone.}
The positive semidefinite cone
$\mathcal K_{\mathrm{psd}} = \{A \in \mathbf{S}^n \, | \, A \succeq 0\}$
has dual cone $\mathcal K_{\mathrm{psd}}$, projection operator 
\[
 \Pi_{\mathrm{psd}}(Z) = \sum_{i=1}^n \max \left\{ 0, \lambda_i \right\} v_i v_i^T
\]
where
$\left\{(\lambda_i, v_i)  \, \middle | \, \lambda_1
\ge \cdots \ge \lambda_i \ge \cdots \ge \lambda_n \right\}$
is the eigendecomposition of $Z$, and 
\[
 D\Pi(Z)[dZ] = V \left( B \circ \left( V^T dZ V \right) \right) V^T,
\]
where $\circ$ denotes the Hadamard (\ie, element-wise) product, $k = \min \left\{ k \, \middle | \, \lambda_k < 0 \right\}$.
The symmetric $B$ is given by
% \[
% B_{ij} = \begin{cases}
%   0 & i \le k, j \le k \\
%   \frac{\max \left\{ \lambda_i, 0 \right\}}
%     {\max \left\{\lambda_i, 0 \right\} - \min \left\{\lambda_j, 0 \right\} }
%     & i > k, j \le k \\
%   \frac{\max \left\{ \lambda_j, 0 \right\}}
%     {\max \left\{\lambda_j, 0 \right\} - \min \left\{\lambda_i, 0 \right\}}
%     & i \le k, j > k \\
%   1 & i > k, j > k.
% \end{cases}
% \]
\[
B_{ij} = \begin{cases}
  0 & i \ge k, j \ge k \\
  \frac{\lambda_i}
    {\lambda_i - \lambda_j }
    & i < k, j \ge k \\
  \frac{\lambda_j}
    {\lambda_j - \lambda_i }
    & i \ge k, j < k \\
  1 & i < k, j < k.
\end{cases}
\]
See~\cite{Busseti2018SolutionRA} for the derivation and proof of $D \Pi_{\mathrm{psd}}$.

\paragraph{Exponential cone.} The exponential cone
\[\mathcal{K}_{\mathrm{exp}} =
\left\{ (x, y, z) \in \reals^3 \, \middle | \, y e^{x/y} \le z, y > 0 \right\}
\cup 
\left\{ (x, 0, z) \in \reals^3 \, \middle | \, x \le 0, z \ge 0 \right\} \]
has dual cone
\[
  \mathcal{K}_{\mathrm{exp}}^* =
  \left\{ (u, v, w) \in \reals^3 \, \middle | \, u < 0, -u e^{v/u} \le ew \right\}
  \cup
  \left\{ (0, v, w) \in \reals^3 \, \middle | \, v \ge 0, w \ge 0 \right\} 
\]
and polar cone $\mathcal{K}_{\mathrm{exp}}^\circ = -\mathcal{K}_{\mathrm{exp}}^*$.
Let $p = (x, y, z)$.
\begin{itemize}
  \item For $p \in \mathcal{K}_{\mathrm{exp}}$, $\Pi_{\mathrm{exp}} (p) = p$.
  \item For $p \in \mathcal{K}_{\mathrm{exp}}^\circ$, $\Pi_{\mathrm{exp}} (p) = 0$.
  \item For $p \not \in \mathcal{K}_{\mathrm{exp}} \cup \mathcal{K}_{\mathrm{exp}}^\circ$
  with $x < 0$ and $y < 0$, $\Pi_{\mathrm{exp}} (p) = (x, 0, \max \left\{ z, 0 \right\})$.
  \item For all other $p$, the projection $\Pi_{\mathrm{exp}} (p)$
    must be found via its definition, \ie, as the (unique) solution to
  \begin{equation}
    \begin{array}{lll}
      &\text{minimize} & \frac{1}{2}\| \hat{p} - p \| \\
      &\text{subject to} & \hat{z} = \hat{y} e^{\hat{x} / \hat{y}}, \, \hat{y} > 0,
    \end{array}
    \label{eq:exp-system}
  \end{equation}
  where $\hat{p} = (\hat{x}, \hat{y}, \hat{z})$ is the optimization variable.
  See~\cite{Friberg2023ProjectionOT} for
  a fast and numerically robust univariate root-finding algorithm that can be used to compute
  the projection.
\end{itemize}
In the following cases, the projection operator is differentiable at $p$.
\begin{itemize}
  \item For $p \in \intr \mathcal{K}_{\mathrm{exp}}$,
  $D \Pi_{\mathrm{exp}}(p) = D \Pi_{\mathrm{exp}}(p)^T = I$.
  \item For $p \in \intr \mathcal{K}_{\mathrm{exp}}^\circ$,
  $D \Pi_{\mathrm{exp}}(p) = D \Pi_{\mathrm{exp}}(p)^T = 0$.
  \item For $p \not \in \mathcal{K}_{\mathrm{exp}} \cup \mathcal{K}_{\mathrm{exp}}^\circ$
  with $x < 0$, $y < 0$, and $z \not = 0$,
  $D\Pi_{\mathrm{exp}}(p) = D\Pi_{\mathrm{exp}}(p)^T =
  \diag (1, 0, \bm{1} \left\{ z > 0 \right\})
  $ where for $\alpha \in \mathbf{R}$, $\bm{1} \left\{ \alpha > 0 \right\}$ is $1$
  if $\alpha > 0$ else is $0$.
  \item For $p \in \intr \left( \reals^3 \setminus \left( \mathcal{K}_{\mathrm{exp}} \cup \mathcal{K}_{\mathrm{exp}}^* 
  \cup \left(\reals_{-} \times \reals_{-} \times \reals \right) \right) \right)$, $D\Pi_{\mathrm{exp}}(p) = (J^{-1})_{1:3, 1:3}$
  where
  \[
    J = 
    \begin{bmatrix}
      1 + \frac{\mu^\star e^{x^\star / y^\star} }{y^\star} & - \frac{\mu^\star x^\star e^{x^\star / y^\star} }{(y^\star)^2}
        & 0 & e^{x^\star / y^\star} \\
      - \frac{\mu^\star x^\star e^{x^\star / y^\star} }{(y^\star)^2} & 1 + - \frac{\mu^\star (x^\star)^2 e^{x^\star / y^\star} }{(y^\star)^3}
        & 0 & (1 - x^\star / y^\star) e^{x^\star / y^\star} \\
      0 & 0 & 1 & -1 \\
      e^{x^\star / y^\star} & (1 - x^\star / y^\star)e^{x^\star / y^\star} & -1 & 0
    \end{bmatrix}.
  \]
  In this Jacobian, $(x^\star, y^\star, z^\star)$ is the solution to~\eqref{eq:exp-system} and $\mu^\star \in \reals$
  is the solution to the dual problem. See~\cite{Ali2017ASN} for the derivation and proof of
  $J$.

\end{itemize}

\paragraph{Dual exponential cone.}
The dual exponential cone is given above.
Via the Moreau decomposition,
its projection operator is $\Pi_{\mathrm{exp}*}(z) =z + \Pi_{\mathrm{exp}}(-z)$
with derivative $D\Pi_{\mathrm{exp}*}(z)[dz] = dz - D\Pi_{\mathrm{exp}}(-z)[dz]$

\paragraph{Power cone.} 
Our power cone and dual power cone results are based on \cite{Hien2015DifferentialPO}.
The (3D) power cone
\[\mathcal{K}_{\mathrm{pow}, \alpha} 
= \left\{ (x, y, z) \in \reals^3 \, \middle | \, x^{\alpha}y^{1-\alpha} \ge \left| z \right|, x \ge 0, y \ge 0  \right\}\]
has dual cone
\[\mathcal{K}_{\mathrm{pow}, \alpha}^{*} =
\left\{ (u, v, w) \in \reals^3 \, \middle | \, \left( \frac{u}{\alpha} \right)^{\alpha}
\left(\frac{v}{1-\alpha} \right)^{1 - \alpha} \ge \left| w \right|, u \ge 0, v \ge 0 \right\}\]
and polar cone $\mathcal{K}_{\mathrm{pow}, \alpha}^{\circ} = - \mathcal{K}_{\mathrm{pow}, \alpha}^{*}$.
%  Its projection operator has no closed-form expression.
Let $p = (x, y, z)$. 
\begin{itemize}
  \item For $p \in \mathcal{K}_{\mathrm{pow}, \alpha}$,
    $\Pi_{\mathrm{pow}, \alpha} (p) = p$.
  \item For $p \in \mathcal{K}_{\mathrm{pow}, \alpha}^\circ$,
    $\Pi_{\mathrm{pow}} (p) = 0$.
  \item For $p \not \in \mathcal{K}_{\text{pow}, \alpha} \cup \mathcal{K}_{\text{pow}, \alpha}^\circ$
    and $z = 0$, $\Pi_{\mathrm{pow}} (p) = (\max \left\{x, 0 \right\}, \max \left\{y, 0 \right\}, 0)$.
  \item For $p \not \in \mathcal{K}_{\text{pow}, \alpha} \cup \mathcal{K}_{\text{pow}, \alpha}^\circ$
    and $z \not = 0$, $\Pi_{\mathrm{pow}} (p) = (f_x, f_y, \textbf{sign}(z)r)$, where
    \begin{align*}
	f_x(r) = \frac{1}{2} \left( x + \sqrt{x^2 + 4 \alpha r (\left| z \right| - r)} \right), \quad
	f_y(r) = \frac{1}{2} \left(y + \sqrt{y^2 + 4 (1 - \alpha) r (\left| z \right| - r)} \right),
    \end{align*}
    and $r$ is the (unique) solution to the (nonconvex) problem
    \begin{equation}
      \begin{array}{lll}
        &\text{find} & r \\
	&\text{subject to} & r = f_x(r)^\alpha f_y(r)^{1 - \alpha} \\
        & & 0 < r < \left| z \right|
      \end{array}
      \label{eq:pow-system}
    \end{equation}
    defined in~\cite[(5)]{Hien2015DifferentialPO}.
\end{itemize}
In the following cases, the projection operator is continuously differentiable at $p$.
\begin{itemize}
\item For $p \in \intr \mathcal{K}_{\mathrm{pow}, \alpha}$,
  $D\Pi_{\mathrm{pow}}(p) = D \Pi_{\mathrm{pow}} (p)^T = I$. 
\item For $p \in \intr \mathcal{K}_{\mathrm{pow}, \alpha}^{\circ}$,
  $D\Pi_{\mathrm{pow}}(p) = D \Pi_{\mathrm{pow}} (p)^T = 0 \in \mathbf{R}^{3 \times 3}$. 
\item \cite[Theorem 3.1]{Hien2015DifferentialPO}
For $p \not \in \mathcal{K}_{\mathrm{pow}, \alpha} \cup \mathcal{K}_{\mathrm{pow}, \alpha}^\circ$
  and $z \not = 0$,
\begin{equation*}
  D\Pi_{\mathrm{pow}}(p) = \begin{bmatrix}
    \frac{1}{2} + \frac{x}{2g_x} + \frac{\alpha^{2} (\left| z \right| - 2r)rL}{g_x^{2}}
    & \frac{(\alpha - \alpha^{2}) (\left| z \right| - 2r)rL}{g_y g_x}
    &\textbf{sign}(z) \frac{\alpha r L}{g_x} \\
    \frac{(\alpha - \alpha^{2}) (\left| z \right| - 2r)rL}{g_x g_y}
    & \frac{1}{2} + \frac{y}{2g_y} + \frac{(1 - \alpha)^2 (\left| z \right| - 2r)rL}{g_y^{2}}
    & \textbf{sign}(z) \frac{(1 - \alpha) r L}{g_y} \\
    \textbf{sign}(z) \frac{\alpha r L}{g_x} & \textbf{sign}(z) \frac{(1 - \alpha) r L}{g_y}
    & \frac{r}{\left| z \right|} + \frac{r}{\left| z \right|}TL
  \end{bmatrix},
  \end{equation*}
  where $g_x = 2f_x - x, g_y = 2f_y - y$,
  \[
    L = \frac{2(\left| z \right| - r)}{\left| z \right| + (\left| z \right| - 2r)(\frac{\alpha x}{g_x}
    + \frac{(1-\alpha)y}{g_y})},
    \quad T = -\left(\frac{\alpha x}{g_x} + \frac{(1 - \alpha)y}{g_y} \right),
  \]
and $r$ is the solution of~\eqref{eq:pow-system}.
 \item For $p \not \in \mathcal{K}_{\mathrm{pow}, \alpha} \cup \mathcal{K}_{\mathrm{pow},\alpha}^\circ$,
 $z = 0$, and $x, y \neq 0$,
  \begin{equation*}
      D\Pi_{\mathrm{pow}} (p) = \begin{bmatrix}
      \bm{1} \left\{ x > 0 \right\} & 0 & 0 \\
      0 & \bm{1} \left\{ y > 0 \right\} & 0 \\
      0 & 0 & d
    \end{bmatrix}.
  \end{equation*}
  The component $d$ is defined as
  \[
    d = \begin{cases}
      1 & x > 0, y < 0, \alpha > 1/2, \text{ or } y > 0, x < 0, \alpha < 1/2 \\
      0 &x > 0, y < 0, \alpha < 1/2, \text{ or } y > 0, x < 0, \alpha > 1/2 \\
      d_x & x > 0, y < 0, \alpha = 1/2 \\
      d_y & x < 0, y > 0, \alpha = 1/2, \\
    \end{cases}
  \]
    where
  \[
    d_x = \frac{x}{2 \left| y \right| + x} \quad \text{and}
    \quad d_y = \frac{y}{2 \left| x \right| + y}.
  \]
\end{itemize}

\paragraph{Dual power cone.}
The dual power cone is given above.
Via the Moreau decomposition,
its projection operator is $\Pi_{\mathrm{pow}*}(z) =z + \Pi_{\mathrm{pow}}(-z)$
with derivative $D\Pi_{\mathrm{pow}*}(z)[dz] = dz - D\Pi_{\mathrm{pow}}(-z)[dz]$

\section*{Acknowledgements}
Parth Nobel was supported in part by the National Science Foundation Graduate Research Fellowship
Program under Grant No. DGE-1656518. Any opinions, findings, and conclusions or recommendations
expressed in this material are those of the author(s) and do not necessarily reflect
the views of the National Science Foundation. The work of Stephen Boyd was supported by
ACCESS (AI Chip Center for Emerging Smart Systems), sponsored by InnoHK funding, Hong
Kong SAR.
\bibliographystyle{abbrv}
\bibliography{./refs}

\end{document}